\documentclass[11pt]{article}
\usepackage[latin1]{inputenc}
\usepackage{graphics}
\usepackage[all]{xy}
\usepackage{amssymb}
\usepackage{amsfonts}
\usepackage[thmmarks]{ntheorem}

\def\RR{{\mathbb R}}

\def\NN{{\mathbb N}}

\def\bea{\begin{eqnarray}}
\def\eea{\end{eqnarray}}

\def\be{\begin{equation}}
\def\ee{\end{equation}}

\newtheorem{theorem}{Theorem}[section]

\theoremstyle{plain}
\theorembodyfont{\normalfont}

\newtheorem{lemma}[theorem]{Lemma}
\theoremstyle{nonumberplain}
\theoremseparator{:}
\theoremsymbol{$\P$}
\newtheorem{demo}{Proof}
\theoremsymbol{$\clubsuit$}

\begin{document}
\author{Fabien Besnard}
\title{An approximation theorem for non-decreasing functions on compact posets}
\maketitle
\abstract{In this short note we prove a theorem of the Stone-Weierstrass sort for subsets of the cone of non-decreasing continuous functions on compact partially ordered sets.}
\section{Introduction}
The classic book \cite{nachbin} contains a theorem which states that given a compact set $M$ and a separating  semi-vector lattice $S$ of continuous real-valued functions on $M$ which contains the constants, there is one and only one way of making $M$ a compact ordered space so that $S$ becomes the set of all non-decreasing continuous real-valued functions on $M$. This theorem has been used in \cite{besnard} to give a putative definition of noncommutative compact ordered sets. However, infimum and supremum turned out to be quite difficult to handle in the noncommutative setting. A different kind of density theorem was thus needed. Since a ``continuous non-decreasing'' functional calculus was available in the noncommutative context, it was natural to look for a theorem which would replace stability under infimum and supremum with stability under continuous non-decreasing functions.

Let us introduce some vocabulary in order to be more precise. Let $S$ be a subset of the set ${\cal C}(M,\RR)$ of continuous real-valued functions on some space $M$, and $H : \RR\rightarrow \RR$ be a function. We will say that $H$ \emph{operates} on $S$ is $H\circ f\in S$ for each $f\in S$.

As remarked in \cite{briem}, the real version of the classical Stone-Weierstrass theorem can be rephrased in terms of operating functions. 

\begin{theorem}(Stone-Weierstrass) Let $S$ be a non-empty subset of ${\cal C}(X,\RR)$, with $X$ a compact Hausdorff space. If  
\begin{enumerate}
\item $S$ is stable by sum,
\item the affine functions from $\RR$ to $\RR$ operate on $S$,
\item the function $t\mapsto t^2$ operates on $S$,
\end{enumerate}
then $S$ is dense in ${\cal C}(X,\RR)$ for the uniform norm.
\end{theorem}

The second hypothesis is a way to say that $S$ is a cone (hence a vector space thanks to first hypothesis) which contains the constant functions.

In fact it is proved in \cite{dlk} that one can replace $t\mapsto t^2$ in the third hypothesis by any continuous non-affine function.

It is a theorem of this kind that we prove in this note, but in the same category (compact ordered sets and non-decreasing continuous functions) as the theorem of Nachbin stated above.

\section{Preliminaries}
Let $M$ be a topological set equipped with a partial order $\preceq $. We let $I(M)$ denote the set of all continuous non-decreasing functions from $M$ to $\RR$, where $\RR$ has the natural topology and the natural ordering, which we write $\le$, as usual. The elements of $I(M)$ are sometimes called continuous isotonies.

Let $S$ be a subset of $I(M)$. We define the relation $\preceq_S$ by

\be
x\preceq _S y\Longleftrightarrow \forall f\in S, f(x)\le f(y)
\ee

It is obvious that $\preceq _S$ is a preorder, which we call \emph{the preorder generated by $S$}. This preorder will be a partial order relation if, and only if, $S$ separates the points of $M$.

We say that \emph{$S$ generates $\preceq$} iff $\preceq _S=\preceq$. This is the case if, and only if, $S$ satisfies 

\be
\forall a,b\in M,\ a\not\preceq b\Longrightarrow \exists f\in S,\ f(a)> f(b)
\ee

Since $a\not=b\Rightarrow a\not\preceq b$ or $b\not\preceq a$, we see that if $S$ generates $\preceq$, it necessarily separates the points of $M$.

Note that it is not guaranteed that for any   poset there exists such an $S$ generating the order. When there is one, then $I(M)$ itself will also generate the order. Posets with the property that $I(M)$ generates the order are called \emph{completely separated ordered sets}. When $M$ is compact and Hausdorff, complete separation is equivalent to the relation $\preceq$ being closed in $M\times M$ (see \cite{nachbin}, p. 114).

%In what follows, we suppose that $M$ is a compact Hausdorff completely separated ordered set.

Let $A$ be a set of functions from $\RR$ to $\RR$. We will say that $A$ operates on $S$ iff 

\be
\forall H\in A, \forall f\in S,\ H\circ f\in S
\ee

%\emph{In what follows we take for $A$ the set of continuous non-decreasing piecewise linear functions.}
\section{Statement and proof of the theorem}

\begin{theorem}
Let $(M,\preceq )$ be a compact Hausdorff partially ordered set. Let $A$ be the set of continuous non-decreasing piecewise linear functions from $\RR$ to $\RR$. Let $S$ be a non empty subset of $I(M)$. If
\begin{enumerate}
\item $S$ is stable by sum,
\item\label{h2} $A$ operates on $S$, 
\item\label{h3} $S$ generates $\preceq$.
\end{enumerate}
then $S$ is dense in $I(M)$ for the uniform norm. 
\end{theorem}

Before proving the theorem, a few comments are in order.

\begin{itemize}
\item First of all, the theorem is true but empty if $M$ is not completely separated, since no $S$ can satisfy the hypotheses in this case.
%\item We also have to be completely clear about what we mean by ``stable under composition with non-decreasing continuous real functions''. It means that for any $f\in S$ and any $H : \RR\rightarrow \RR$, continuous and non-decreasing, then $H\circ f\in S$. It is maybe useful to note that it is equivalent to require that for any $f\in S$ and any continuous non-decreasing $H$ from $f(M)$ to $\RR$ one has $H\circ f\in S$. Indeed, any continuous non-decreasing function from the compact set $f(M)$ can be extended to a continuous non-decreasing function defined on the whole of $\RR$ (to see this one can use the fact that $\RR$ is normally ordered, see theorem 7 p 103 in \cite{nachbin}).
\item The hypothesis that $S$ is not empty is redundant if $M$ has at least two elements, by \ref{h3}.
\item Finally, let us remark that \ref{h2} entails that $S$ is in fact a convex cone which contains the constant functions.
\end{itemize}
%(since $f(M)$ in included in a toposet $[a,b]$ which has the Tietze-Nachbin extension property, and we can then extend to $\RR$ by the constant $f(b)$ to the right and $f(a)$ to the left).

%To prove the theorem we can as well suppose that $S$ is closed since the closure $\bar S$ will also satisfy the hypotheses if $S$ does.

To prove the theorem we  need two lemmas.

\begin{lemma}\label{premierlemme} Let $x,y\in M$ be such that $y\not\preceq  x$. Then $\exists f_{x,y}\in S$ such that $0\le f_{x,y}\le 1$, $f_{x,y}(x)=0$ and $f_{x,y}(y)=1$.
\end{lemma}
\begin{demo}
Since $S$ generates $\preceq$, there exists $f\in S$ such that $f(x)<f(y)$. Let $H\in A$ be such that $H(t)=0$ for $t\le f(x)$, $H$ is affine on the segment $[f(x),f(y)]$, and $H(t)=1$ for $t\ge f(y)$. Then  $f_{x,y}:=H\circ f$ meets the requirements of the lemma.
\end{demo}

\begin{lemma}\label{secondlemme} Let $K,L$ be two compact subsets of $M$ such that $\forall x\in K,\forall y\in L$, $y\not\preceq x$. Then $\exists f_{K,L}\in S$ such that $0\le f_{K,L}\le 1$, $f=0$ on $K$ and $f=1$ on $L$.
\end{lemma}
\begin{demo} For all $x\in K$ and $y\in L$, we find a $f_{x,y}\in S$ as in lemma \ref{premierlemme}. We fix a $y\in L$ and let $x$ vary in $K$. Since $f_{x,y}$ is continuous, there exists an open neighbourhoud $V_x$ of $x$ such that $f_{x,y}(V_x)\subset [0;1/4[$. By compacity of $K$, there exists $V_1,\ldots,V_k$ corresponding to $x_1,\ldots,x_k$ such that $K\subset V_1\cup\ldots\cup V_k$.

Now we define $g_y:={1\over k}\sum_if_{x_i,y}$. We have $g_y\in S$ since $S$ is a convex cone (see the last remark  below the theorem). It is clear that $g_y(y)=1$ and that for all $x\in K$, $0\le g_y(x)\le {1\over k}(k-1+1/4)=1-{3\over 4k}<1$. We then choose $H\in A$  such that $H(t)=0$ for  $t\le 1-{3\over 4k}$ and $H(t)=1$ for $t\ge 1$. We set $f_{K,y}:=H\circ g_y$. We thus have $f_{K,y}\in S$,  $f_{K,y}=0$ on $K$ and $f_{K,y}(y)=1$.

Using the continuity of $f_{K,y}$, we find an open neighbourhood $W_y$ of $y$ such that $f_{K,y}(W_y)\subset[3/4;1]$. Since we can do this for every $y\in L$, and since $L$ is compact, we can find functions $f_{K,y_j}$, $j=1..l$, and open sets $W_1,\ldots,W_l$ of the above kind such that $L\subset W_1\cup\ldots\cup W_l$. We then define $g={1\over l}\sum_jf_{K,y_j}$. We have $g\in S$,  and $g(K)=\{0\}$. Moreover, for all $z\in L$, $1\ge g(z)\ge {3\over 4l}>0$. We then  choose a function $G\in A$  such that $G(t)=1$ for $t\ge 3/4l$ and $G(t)=0$ for $t\le 0$. Now the function $f_{K,L}:=G\circ g$ has the desired properties.
\end{demo}

We can now prove the theorem.

\begin{demo}
Let $f\in I(M)$. We will show that, for all $n\in\NN^*$ there exists $F\in S$ such that $\|f-F\|_\infty\le{1\over n}$.

If $f$ is constant then the result is obvious. Else, let $m$ be the infimum of $f$ and $M$ be its supremum. Let $\tilde f={1\over M-m}(f-m.1)$. Using the fact that $S$ is a convex cone, we can work with $\tilde f$ instead of $f$. Hence, we can suppose that $f(M)=[0;1]$ without loss of generality. 

We set $K_i=f^{-1}([0;{i\over n}])$, and $L_i=f^{-1}([{i+1\over n};1])$ for each $i\in\{0;\ldots;n-1\}$. Since $f$ is continuous and $M$ is compact, the sets $K_i$ and $L_i$ are both closed, hence compact.

For each $i$ we use lemma \ref{secondlemme} to find $f_i\in S$ such that $f_i(K_i)=\{0\}$ and $f_i(L_i)=\{1\}$.

We then consider the function $F={1\over n}\sum_{i=0}^{n-1}f_i$. We clearly have $F\in S$.

Let $m\in M$. Suppose ${j\over n}<f(m)<{j+1\over n}$ for some $j\in\{0;\ldots;n-1\}$. We thus have $m\in K_i$ for $j< i< n$ and $m\in L_i$ for $i<j$. Hence $F(m)={1\over n}\sum_{i=0}^jf_i(m)={1\over n}(j+f_j(m))\in [{j\over n};{j+1\over n}]$. Thus $|f(m)-F(m)|\le {1\over n}$.

Now suppose $f(m)={j\over n}$, with $j\in\{0;\ldots;n\}$. We have $m\in K_i$ for $i\ge j$ and $m\in L_i$ for $i<j$. Thus $F(m)={1\over n}\sum_{i=0}^{j-1}f_i(m)={j\over n}$. We see that $|f(m)-F(m)|=0$ in this case.

Hence we have shown that $|f(m)-F(m)|\le {1\over n}$ for all $m\in M$, thus proving the theorem.
\end{demo}

To conclude, let us remark that the set $A$ in the theorem can be replaced by any subset of $I(M)$ with the following property : for any two reals $a<b$, there exists $f\in A$ such that $f=0$ on $]-\infty;a]$, and $f=1$ on $[b;+\infty[$. For example one can take $A=I(M)$ itself, or the set of non-decreasing ${\cal C}^\infty$ functions.
% \section{Final Comments}
% As remarked in \cite{briem}, the classical Stone-Weierstrass theorem can be rephrased

% It is obvious from the proof that only the real functions which are continuous non-decreasing and affine by segment are really needed in the second hypothesis of the theorem. It is an exercise to show that we could also do with non-decreasing polynomials. It would be interesting to find other classes of functions that could be used.

\end{document}